\documentclass{article}

\usepackage{arxiv}

\usepackage[utf8]{inputenc} 
\usepackage[T1]{fontenc}    
\usepackage{hyperref}       
\usepackage{url}            
\usepackage{booktabs}       
\usepackage{amsfonts}       
\usepackage{nicefrac}       
\usepackage{microtype}      
\usepackage{lipsum}
\usepackage{cite}
\usepackage{times} 
\usepackage{amsmath} 
\usepackage{amssymb}
\usepackage{graphicx}
\graphicspath{{Figures/}}

\title{Identification and Validation of Virtual Battery Model for Heterogeneous Devices
}

\author{
  Sai Pushpak Nandanoori, Indrasis Chakraborty, Thiagarajan Ramachandran and Soumya Kundu\thanks{The authors would like to thank U.S. Department of Energy for supporting this work under the ENERGISE program (contract DE-AC02-76RL01830). The authors are with the Optimization and Control Group, Pacific Northwest National Laboratory, Richland, WA 99354 USA. Emails: \{saipushpak.n, indrasis.chakraborty, thiagarajan.ramchandran, soumya.kundu\}@pnnl.gov.}
}

\begin{document}
\maketitle

\begin{abstract}
The potential of distributed energy resources in providing grid services can be maximized with the recent advancements in demand side control. Effective utilization of this control strategy requires the knowledge of aggregate flexibility of the distributed energy resources (DERs). Recent works have shown that the aggregate flexibility of DERs can be modeled as a virtual battery (VB) whose state evolution is governed by a first order system including self dissipation. The VB parameters (self dissipation rate, energy capacity) are obtained by solving an optimization problem which minimizes the tracking performance of the ensemble and the proposed first order model. For the identified first order model, time varying power limits are calculated using binary search algorithms. Finally, this proposed framework is demonstrated for different homogeneous and heterogeneous ensembles consisting of air conditioners (ACs) and electric water heaters (EWHs).
\end{abstract}

\section{Introduction}
The last few years have seen a significant increase in integration of renewable energy into the electricity grid, and the intermittent nature of renewable energy causes uncertainty in power generation. Moreover, on the distribution side, there is an increased visibility of thermostatically controllable loads (TCLs) such as ACs and EWHs due to advancements in power electronics, communication capabilities that enable remote monitoring/controlling of TCLs. With increased renewable penetration, these advancements allow TCLs to provide several grid services such as demand response, frequency regulation, frequency response, and tracking regulation \cite{ callaway2009tapping, callaway2011achieving, koch2011modeling, hao2015aggregate, bashash2011modeling, zhang2013aggregated, kundu2011modeling, pushpak2018prioritized}. In this work, we focus on characterizing the aggregate flexibility of the homogeneous and heterogeneous ensemble of ACs and EWHs using a first-order virtual battery (VB) model. This work is a preliminary step toward characterizing heterogeneous energy resources such as ACs, EWHs, Solar PV, electric vehicles and real battery as a single VB model. 

The allowable power flexibility for any device is characterized by a set whose values are allowable deviations from the nominal power. Furthermore, the aggregate flexibility of the ensemble of heterogeneous devices is obtained by taking the Minkowski sum of the flexibility sets corresponding to each device \cite{ kundu2018approximating}. Due to the complexity involved in identifying the exact Minkowski sum of the ensemble, the authors in \cite{kundu2018approximating} have come up with a method that gives an outer estimate of the aggregate flexibility using geometric methods. Further, authors in \cite{indrasis2018geometric} have extended the geometric method by incorporating the uncertainties into characterizing the stochastic aggregate flexibility of the ensemble. Another way of characterizing the aggregate flexibility of the ensemble is by representing it with an equivalent possibly a first-order linear system, VB \cite{hao2015aggregate,hao2017optimal,hughes2016identification}. The state of charge (\textit{soc}) of the VB denotes the aggregate flexibility of the ensemble and is captured simply in terms of battery parameters such as power limits, initial \textit{soc}, self dissipation rate and capacity. 

There exists several characterizations of aggregate flexibility using VB models in the literature \cite{ mathieu2015arbitraging, hao2015aggregate, hughes2016identification, hao2017optimal}. These works include characterizing a VB model for a wide range of systems from small residential TCLs to complex systems such as commercial building heating, ventilating, and air-conditioning (HVAC) loads. The thermal mass in the residential and commercial buildings enables them to store thermal energy and be used to provide grid services such as frequency regulation by varying the device power consumption. Traditionally, such services are provided by varying the generation at the generators. Similar to a real battery, the VB also has self dissipation rate, capacity, and power limits as parameters. The objective is to identify these VB parameters corresponding to the collection of homogeneous or heterogeneous TCL devices. 

Two of the recent and most interesting works in this area are by Hughes, et al. \cite{ hughes2015virtual,  hughes2016identification}. The authors compare the characteristics of a linear VB model with the actual nonlinear ensemble to identify the VB parameters by solving an optimization problem. This is interesting in the sense that they presented a realistic case study to illustrate the potency of the proposed approach. The VB model discussed here is inspired by \cite{hughes2016identification} and follows it closely to characterize a VB model for an ensemble of homogeneous and heterogeneous devices. 

The main contributions of this work can be summarized as follows. We formulated an optimization problem for identifying the device state (ON/OFF) in any given ensemble, such that the ensemble can successfully track the selected regulation signal, while satisfying the device level temperature constraints. Utilizing device level information from the aforementioned optimization problem, the VB parameters are identified by minimizing error in performance between the actual device ensemble and the equivalent VB model. In contrast to the existing works, we propose analytic expressions for the initial estimate of the VB state for different type of ensembles. Binary search algorithms are used to identify time-varying power limits of the identified VB model. Finally, applying the proposed framework, VB system is identified and validated for different homogeneous and heterogeneous ensemble of ACs and EWHs.

\section{Virtual Battery Characterization}
\label{sec:VB_model}
We begin this section with the description of the dynamics of the first order VB system. 
\begin{subequations}
\begin{align}
\dot{x}(t) &= - a x(t) - u(t)\,,\quad x(0) = x_0 \label{eq:VB_model} \\
C_1&\leq x(t) \leq C_2, \label{eq:VB_capacity} \\
P^-(t) &\leq u(t) \leq P^+(t),\label{eq:VB_constraints}
\end{align}\end{subequations}
where $x(t) \in \mathbb{R}$ denotes the state of charge (\textit{soc}) of the VB at time $t$ with the initial \textit{soc} $x_{0}$, $a$ denotes the self dissipation rate, and the lower and upper energy limits of the VB are denoted by $C_1$ and $C_2$, respectively. The regulation signal, $u(t)$ acts as an input to the VB and must always lie within the time-varying power limits $P^{-}(t)$ and $P^+(t)$. This first-order VB model can be applied to characterize the aggregated flexibility of distributed energy resources (DERs) and many building loads \cite{hao2015aggregate,hao2017optimal,hughes2016identification}, where the unknown VB parameters are denoted by vector $\phi = [a, C_1, C_2, x_0, P^-(t) ,P^+(t)]$.

Homogeneous and heterogeneous ensemble of DERs can be represented using this first order VB model. For example, if an electric vehicle or a real battery is considered, then they appear parallel to the VB and an equivalent battery model can be derived. However identifying VB parameters when an electric vehicle is in the ensemble, requires further investigation, due to presence of additional constraints (such as charging constraints), coming from the electric vehicle. The inclusion of solar PV results only in additional power limits as this is a free instantaneous source of energy, and hence does not result in an increased battery capacity. Although this scenario changes if solar PV is connected to a battery storage. The scope of this work extends to ACs, EWHs, and a heterogeneous ensemble of ACs and EWHs. We refer the readers to \cite{indrasis2018transferlearning, pushpak2018prioritized} for the first order hybrid models of ACs and EWHs. In what follows, we describe the method to identify the VB parameters, $\phi$.

\section{Identification of VB Parameters}
\label{sec:VB_model_identification}
The VB parameters are computed in two steps: the first step involves finding the time at which the ensemble fails to track a regulation signal and the second step minimizes the difference in performance between the actual DER ensemble and the equivalent VB model \cite{hughes2016identification}. We begin the discussion on control design problem to track a regulation signal.

\subsection{Tracking regulation signals}
For a given ensemble of ACs/EWHs, the solution to the following optimization problem gives a control law indicating which devices have to be on ON/OFF state such that the regulation signal is tracked at all times and user defined set-points are satisfied. 
\begin{align}
\mbox{minimize}_s \; & \vert \vert T(t+1) - T^{set}\vert \vert^2_2 + \vert \vert T_w(t+1) - T^{set}_w\vert \vert^2_2 \nonumber \\
\mbox{subject to}\; & T^{set} - \delta T/2 \leq T(t+1) \leq T^{set} + \delta T/2 \label{eq:constraint_AC}\\
& T_w^{set} \!-\! \delta T_w/2 \!\leq\! T_w(t+1) \!\leq\! T_w^{set} + \delta T_w/2 \label{eq:constraint_EWH} \\
& \vert u_i(t) - s^{\top} P \vert \leq \varepsilon \label{eq:constraint_power}\\
& s \in \{0,1\} \label{eq:constraint_device_status}
\end{align}
where $\varepsilon > 0$, $u_i(t)$ is the $i^{th}$ regulation signal at time $t$ and $T(t), T_w(t)$ denote the vector of temperatures corresponding to AC and EWH devices at time $t$ respectively. The temperature set points and dead band temperature limits for AC (EWH) devices are given by $T^{set}$ ($T^{set}_w$) and $\delta T (\delta T_w)$. 

The constraints, Eqs. \eqref{eq:constraint_AC} and \eqref{eq:constraint_EWH} refer to the temperature bounds on AC and EWH devices. Notice that this optimization problem is nonconvex due to the binary constraints on the optimization variable $s_i$ (ON/OFF status of each device). We make the following two assumptions on the optimization problem to make it computationally tractable.  
\begin{itemize}
    \item relax the binary constraints on $s_i$ to a convex constraint, $0 \leq s_i \leq 1$ for all $i$, and
    \item adding the additional cost on $s_i$, i.e., $\sum_i s_i (1-s_i)$.
\end{itemize}

The binary relaxation of the optimization variable $s_i$ ensures that the optimization problem is convex. This results in a value of $s_i$ lying in the set $[0,1]$ which essentially implies that the power applied to either AC or EWH device is continuous which is not correct (refer Eqs. (8) and (9) in \cite{indrasis2018transferlearning} or Eqs. (11) and (14) in \cite{pushpak2018prioritized}). Hence, we introduce the additional (nonconvex) cost on the optimization variable, $s_i$ by forcing it to take either $0$ or $1$ value which results in the following control design problem. 
\begin{equation} 
\begin{aligned}
    & \mbox{minimize}_s \qquad   W_1 \vert \vert T(t+1) - T^{set}\vert \vert^2_2  + \textstyle{\sum_i s_i (1-s_i)} + W_2 \vert \vert T_w(t+1) - T^{set}_w\vert \vert^2_2  \\
    & \mbox{subject to} \qquad   0 \leq s \leq 1, \text{ and} \quad  \mbox{Eqs}.\; \eqref{eq:constraint_AC}-\eqref{eq:constraint_power}  
\end{aligned}
\label{eq:optim_problem_convex}
\end{equation}
where $W_1$ and $W_2$ are some positive weights on the respective costs.
\begin{figure} 
\begin{center}
\includegraphics[width = 0.6 \textwidth]{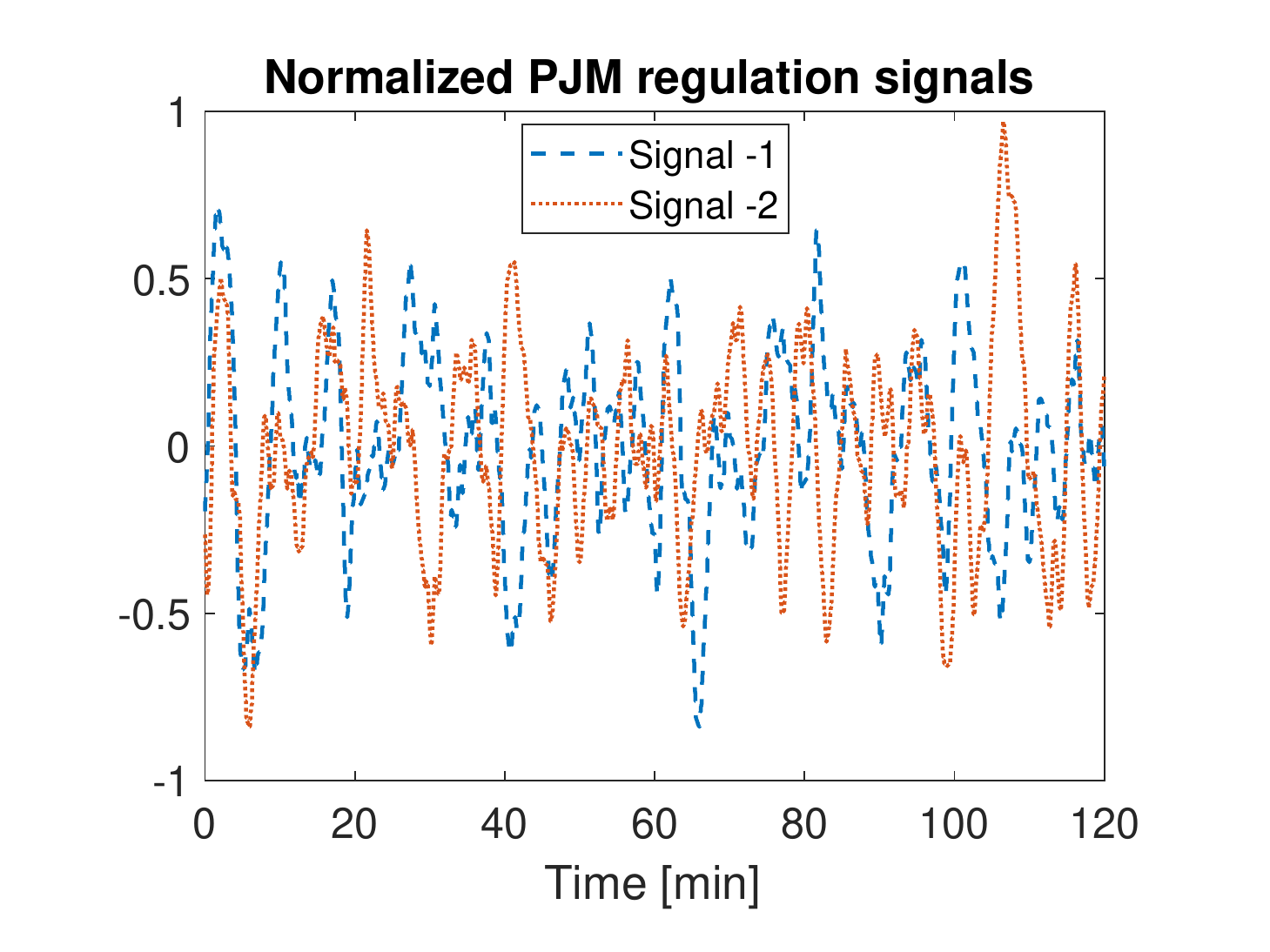}
\caption{Two sample PJM signals}
\label{fig:reg_signals}
\end{center}
\end{figure}
For each regulation signal, the control design problem, Eq. \eqref{eq:optim_problem_convex} is solved to identify whether the ensemble can be able to track the regulation signal at every time step. If at any time point, the solution does not exist then it implies that the ensemble failed to track the regulation signal. Although this process results in increased computation time, it ensures that the VB model thus obtained is accurate. PJM regulation signals \cite{PJM_signals} are considered in this work and sample normalized regulation signals corresponding to a two-hour time window can be seen in Fig. \ref{fig:reg_signals}.  

Any feasible regulation signal corresponding to an ensemble is computed by adding the normalized regulation signals to the baseline power consumption of the ensemble. The baseline power consumption is the actual power consumption of the ensemble without any regulation signal being applied. More discussion on this computation is reserved until Section \ref{sec:simulation}. 
 
The criterion for finding the VB model parameters is a comparison of the violation times of the ensemble and the linear first-order VB model. We recall the notion of violation times \cite{hughes2016identification,hughes2015virtual} in the following subsection. 
 
\subsection{Violation times}
The violation time of the ensemble is defined as the minimum time instant at which the control design problem, Eq. \eqref{eq:optim_problem_convex} fails to find a solution \cite{hughes2016identification}. Similarly, the violation times for the first-order VB model can be computed when constraints given in Eq. \eqref{eq:VB_constraints} are violated. Corresponding to the regulation signal, $u_i(t)$, the violation time function for the ensemble and VB are respectively denoted by $f(u_i(t), \bar \tau)$ and $b(u_i(t), \phi,\bar \tau)$. 

In practical computation of these violation times, a time horizon is chosen over which the regulation signal is considered. Therefore when the ensemble succeeds in tracking the regulation signal, the violation time is considered to be the maximum time in which the regulation signal exists. We refer the readers to \cite{hughes2016identification,hughes2015virtual} for a detailed discussion on computation of violation times. Next, the discussion on VB parameter identification follows.

\subsection{VB parameter identification}
The time varying power limits of the VB are identified by extending the binary search algorithm proposed in \cite{hughes2015virtual}. The algorithm to compute $P^+(t)$ is shown below. 
\begin{center}
\begin{tabular}{p {0.8 \textwidth}}
\hline 
Algorithm: Finding upper power limit, $P^+(t)$ at each $t$ \\
\hline
Initialize $\alpha = 0$ \hfill \textit{Lower bound} \\
Initialize $\beta = 1$  \hfill \textit{Upper bound} \\
\textbf{while} the ensemble tracks $\hat u(t) := \beta + P_{base}(t)$ \hfill \textit{No instant violation} \\
\quad $\alpha := \beta$, $\beta := 2 \times \beta$ \hfill \textit{Update $\alpha$ and $\beta$} \\
end \textbf{while} \\
\textbf{while} $(\beta - \alpha) > \epsilon$ \hfill \textit{Testing a new bound} \\
\quad $\gamma := \frac{\alpha + \beta }{2}$ \\
\quad \textbf{if} ensemble tracks $\hat u(t) := \gamma + P_{base}(t)$, then $\alpha := \gamma$ \\
\quad \textbf{else} $\beta := \gamma$ \\
\quad end \textbf{if} \\
end \textbf{while} \\
return $P^+(t) := \alpha$. \\
\hline
\end{tabular}
\end{center}
Similarly, the algorithm to find $P^-(t)$ can be formulated by changing  $\beta$ to $-\beta$ and $\gamma$ to $-\gamma$ in the above algorithm. The initial condition of the VB system corresponding to the AC ensemble is denoted by $x_0 = x(t)|_{t=0}$ where 
\begin{align}
x(t) = \sum_i \frac { T_i(t) - (T_i^{set}-\delta T/2)}{\eta_i/C_i }
\label{eq:VB_AC_ic0}
\end{align}
Similarly, the initial condition to the first-order VB model corresponding to an EWH ensemble is given by $x_{w0} = x_w(t)|_{t=0}$ where 
\begin{align}
x_w(t) = \sum_i \frac { T_{w_i}^{set}+\delta T_w/2 - T_{w_i(t)} }{1/C_{w_i} }
\label{eq:VB_WH_ic0}
\end{align}
Finally in the case of VB for a heterogeneous ensemble (ACs and EWHs), the initial condition for the VB can be obtained as a sum of initial conditions defined in Eqs. \eqref{eq:VB_AC_ic0} and \eqref{eq:VB_WH_ic0}. \looseness=-1

The self dissipation rate and battery capacity bounds of the VB are obtained as a solution to the following optimization problem \cite{hughes2016identification}. \looseness=-1
\begin{equation}
\begin{aligned}
& \text{minimize}_{a, C_1, C_2} \quad  ||\, \log|B(a, C_1, C_2, x_0) - F|\, ||_2 \\
& \text{subject to} \quad \quad  B(a,C_1, C_2,x_0) \leq F
\end{aligned}  
\label{eq:VB_optm_problem}
\end{equation}
where the vectors $B$ and $F$ are given by 
\[B =  [ b(u_1(t),\phi,\bar \tau)  \dots  b(u_n(t),\phi,\bar \tau) ]^{\top},\] 
\[F =  [ f(u_1(t),\bar \tau)  \dots  f(u_n(t),\bar \tau) ]^{\top} 
\]
and $n$ indicates the number of regulation signals. The cost in this optimization problem always tries to minimize the difference between the violation times of the linear first-order VB system and the actual ensemble. The motivation behind the optimization problem, \eqref{eq:VB_optm_problem}, is that if the first-order VB model does not cause any violations, then the ensemble also does not cause any violations. 

Note that, in \cite{hughes2016identification}, the initial condition $x_0$ of the VB system is identified as a solution to the optimization problem given in \eqref{eq:VB_optm_problem}. By making the initial \textit{soc} of the VB an optimization variable, the VB parameters thus obtained are conservative. In this work, we overcome this conservatism by giving an analytical expression for computing the initial \textit{soc} for ACs and EWHs. 
In the following section, we consider an ensemble of AC and EWH devices and identify their corresponding VB parameters by applying the proposed framework. 

\section{Simulation Study}
\label{sec:simulation}
In this section, we present the VB parameters for different ensembles and compare the values thus obtained with existing VB parameter computation methods. The first step in computing the VB parameters of an ensemble is to compute the actual violation times by solving the control design optimization problem (shown in Eq. \eqref{eq:optim_problem_convex}) with respect to a set of regulation signals. The regulation signals are computed by combining the baseline power consumption and normalized PJM regulation signals. Discussion on the normalized PJM regulation signals is given in Section \ref{sec:VB_model_identification}. 

\begin{figure}
\begin{center}
\includegraphics[width = 0.65 \linewidth ]{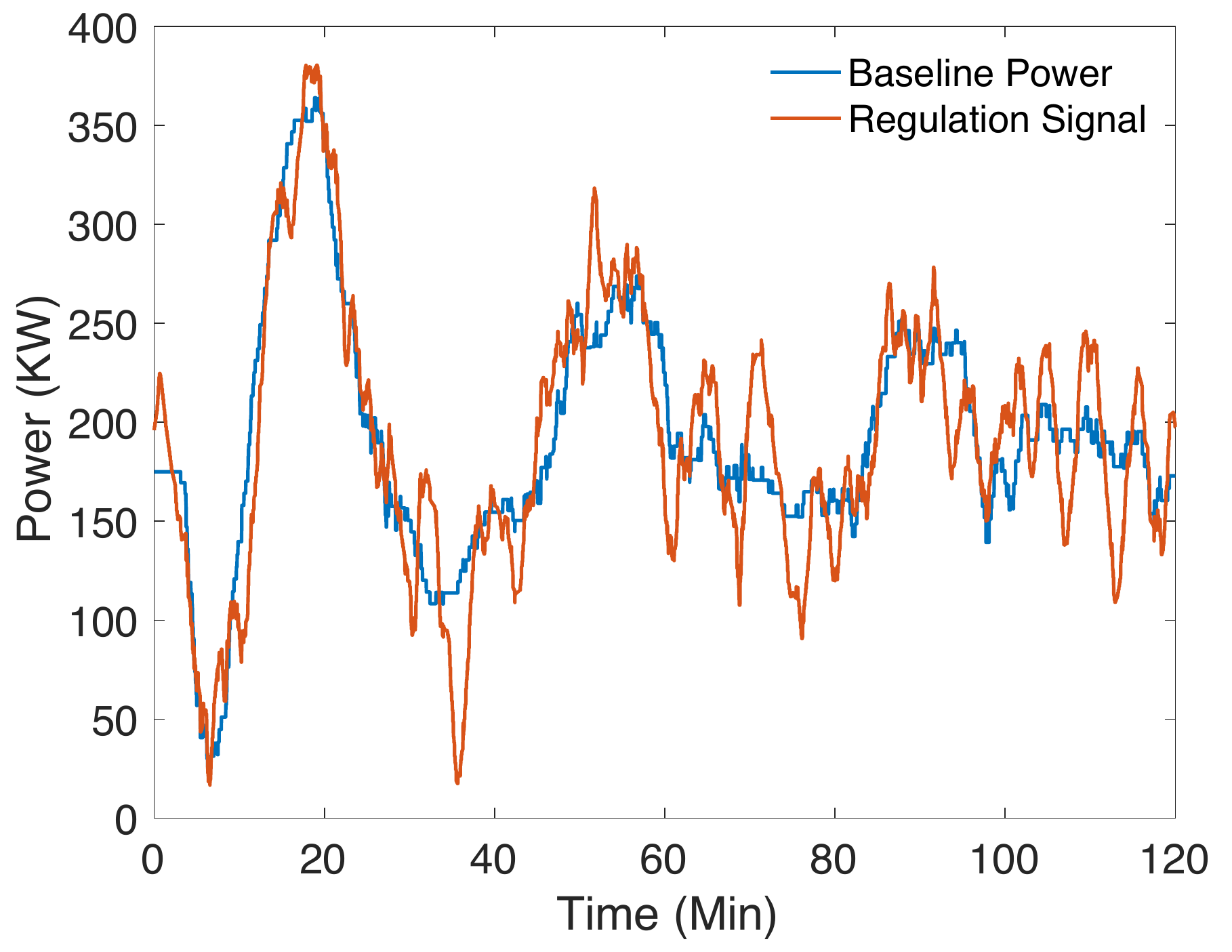}
\caption{Baseline power of the ensemble and a regulation signal}
\label{fig:reg_signal_and_baseline_power}
\end{center}
\end{figure}
The baseline power consumption for the ensemble of AC devices is computed by allowing the ensemble to evolve in the absence of a regulation signal. In this work, a time window of 2 hours is considered and hence the baseline power is computed for 2 hour time horizon. The mathematical representation for the computation of regulation signal is given by 
\begin{align*}
P_{reg} = P_{base} + \tilde P_{reg}
\end{align*}
where $\tilde P_{reg} = \gamma P_{base}^{mean} \hat P_{reg}^{PJM}$; $\gamma$ is the regulation factor; $ P_{base}^{mean} $ is the average of the baseline power; and $\hat P_{reg}^{PJM}$ is the normalized PJM regulation signal. This computation ensures that the regulation signal obtained in this way is a feasible input for the ensemble. The regulation signal and baseline power corresponding to an ensemble of 100 AC devices can be seen in Fig. \ref{fig:reg_signal_and_baseline_power}. The 2-hour time window is arbitrary and the resolution of the regulation signal is chosen to be $1$ second. 

\subsection{VB Model: Homogeneous Devices}
Here we discuss the computation of VB parameters for two types of homogeneous devices, ACs and EWHs. The AC ensemble consists of $100$ devices and EWH ensemble consists of $120$ devices. The temperature evolution of AC and EWH devices is given in Eqs. (8) and (9)  \cite{indrasis2018transferlearning} respectively. 
The parameters such as set-point temperature, outside air temperature, thermal resistance and thermal capacitance for the $100$ AC devices is chosen such that no two devices have the same parameters. Similarly, the device parameters for all the EWHs are chosen. 
A medium water flow profile is considered for all the EWHs. The baseline power corresponding to $100$ AC devices is shown in Fig. \ref{fig:reg_signal_and_baseline_power}.

\begin{figure}
    \centering
    \includegraphics[width = 0.65 \linewidth]{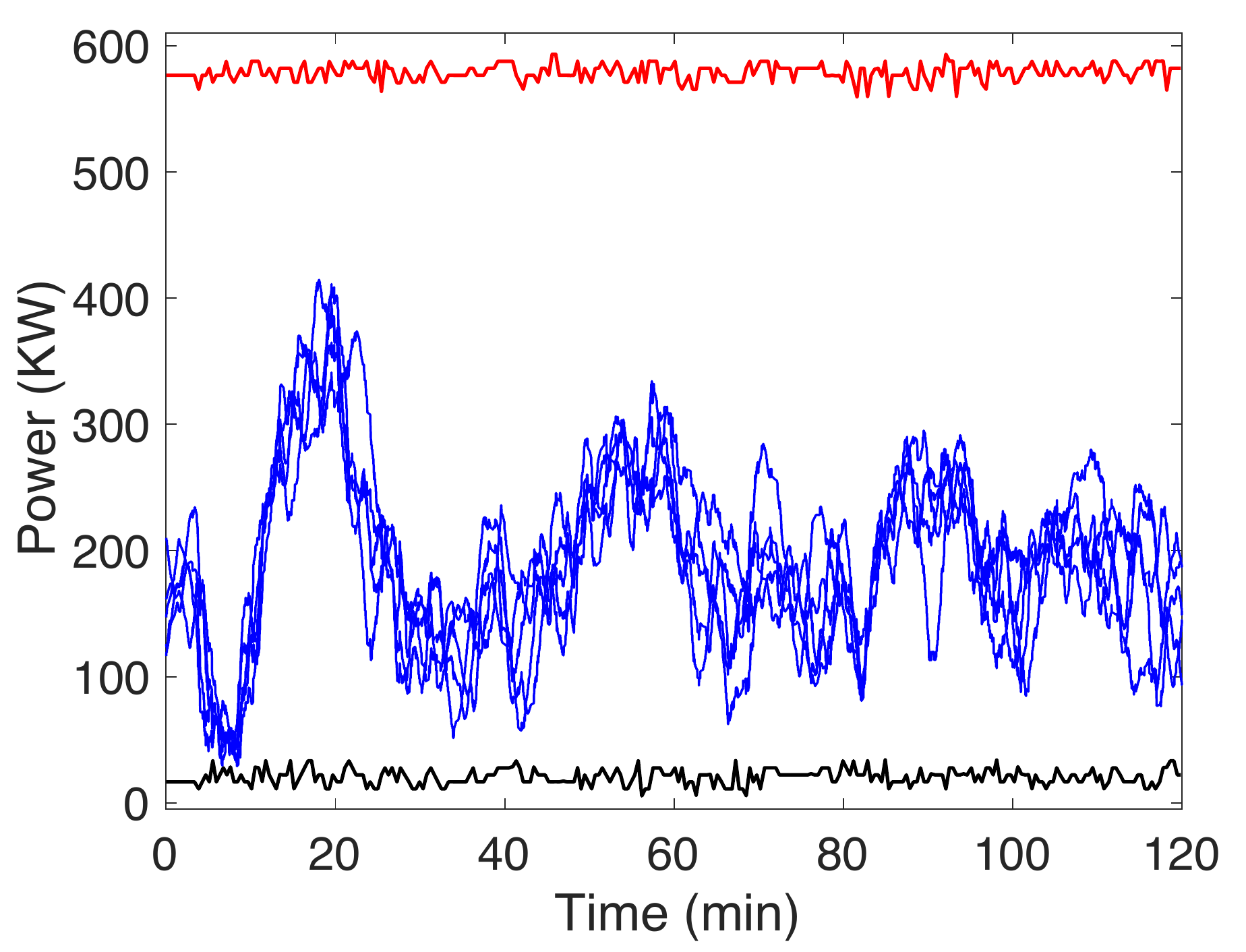}
    \caption{Time varying power limits for 100 AC devices and PJM regulation signals}
    \label{fig:power_limits}
\end{figure}
We first identify the time-varying power limits for each ensemble. The upper and lower power limits ($P^{+}(t)$ and $P^{-}(t)$) for AC and EWH ensembles are computed by applying the binary search algorithm discussed in Section \ref{sec:VB_model_identification}. Each regulation signal must satisfy the power limits before asking the ensemble to track it.

Here, we consider $200$ regulation signals for each ensemble and the ones violating the power limits are discarded. The remaining set of regulation signals are then used to identify the violation times. Time-varying power limits corresponding to the AC ensemble and some regulation signals  satisfying the power limits can be seen in Fig. \ref{fig:power_limits}. 
Next, the control design optimization problem given in Eq. \eqref{eq:optim_problem_convex} is solved in Julia by using the Ipopt solver to obtain the violation time wrt each regulation signal. The solution to the optimization problem (Eq. \eqref{eq:optim_problem_convex}) identifies which AC device has to be in ON or OFF state such that the regulation signal is tracked by the ensemble. As a result, the violation times are computed for AC as well as EWH ensembles.

The violation times for AC and EWH ensemble are now used to solve for the VB parameters such as self dissipation rate and capacity (see optimization problem, Eq. \eqref{eq:VB_optm_problem}). It can be seen from the optimization problem given in Eq. \eqref{eq:VB_optm_problem} that the optimization variable, $s_i$ which decides the ON or OFF status of any device is allowed to take values between 0 and 1. However, the additional cost on the optimization variable, $s_i$ ensures that the resultant values are close to either $0$ or $1$ and thereby trying to reduce the error due to convex relaxation.
\begin{figure}
    \centering
    \includegraphics[width = 0.65 \linewidth]{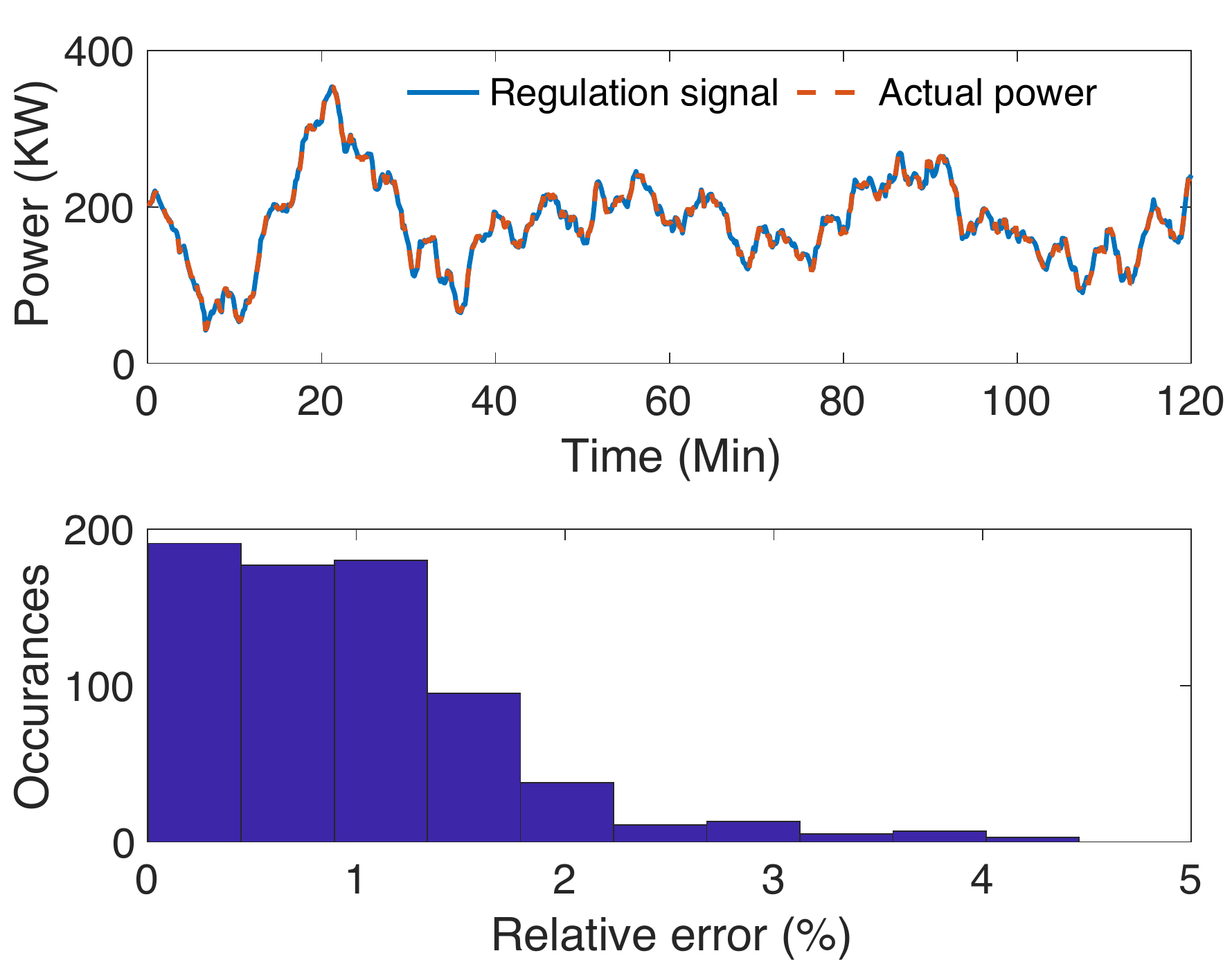}
    \caption{(a) Regulation signal and the AC ensemble power in a 2 hour time window (b) Relative error histogram}
    \label{fig:tracking_and_error}
\end{figure}
Therefore, to identify the (binary) status of each device, the solution to the optimization problem Eq. \eqref{eq:optim_problem_convex}, $s_i$ is rounded off to its closest integer. Thus the total power of the ensemble can be computed and is shown along with the regulation signal in Fig. \ref{fig:tracking_and_error}(a).

Furthermore the distribution of percentage relative error between the regulation signal and the ensemble power is shown in Fig. \ref{fig:tracking_and_error}(b). 
\begin{figure}
\begin{center}
    \includegraphics[width = 0.65 \linewidth]{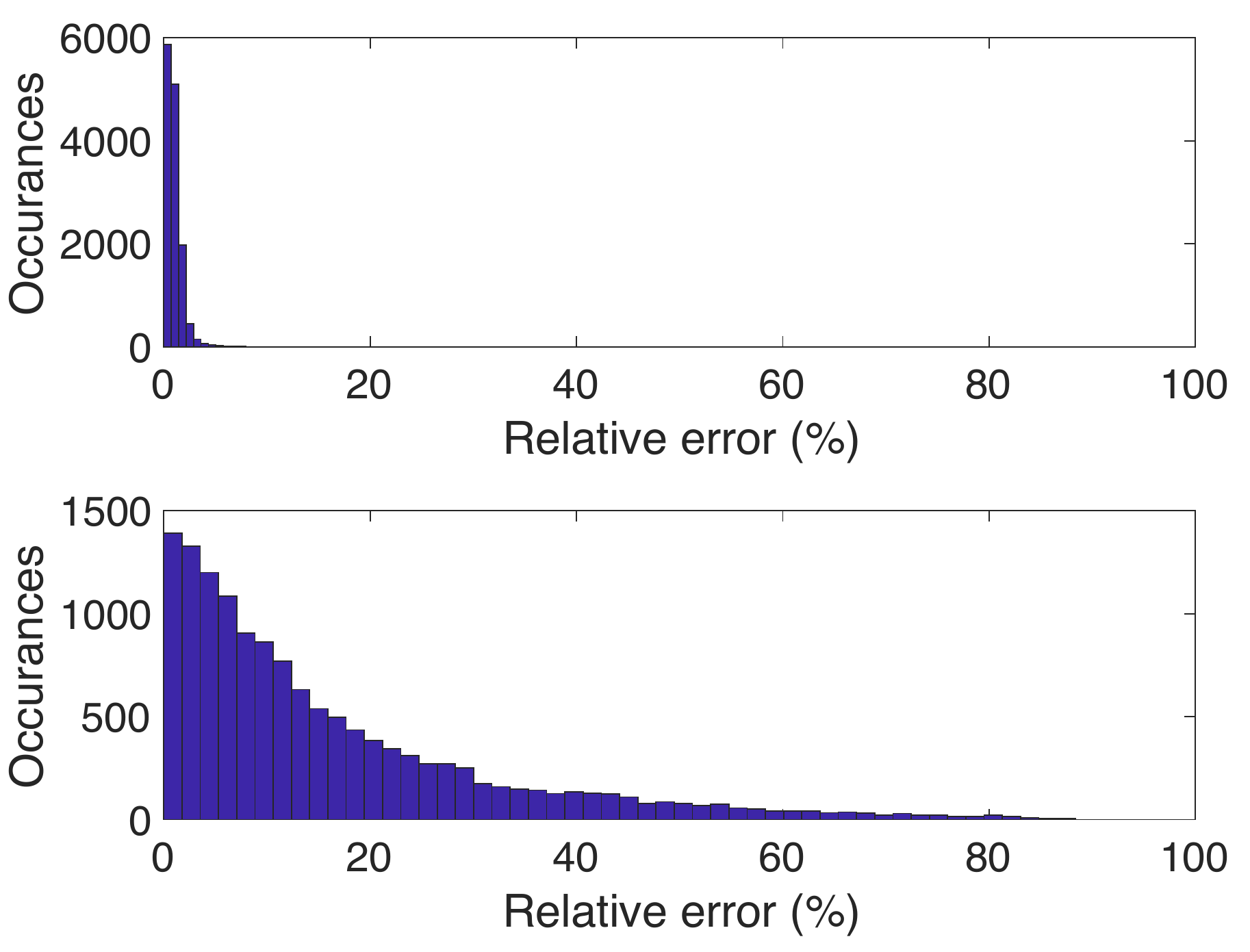}
    \caption{Histogram of percentage error between the ensemble power and various regulation signals (a) With additional cost on $s$ (b) No cost on $s$.}
    \label{fig:full_error_hist}
    \end{center}
\end{figure}
From Fig. \ref{fig:tracking_and_error}(b), it can be seen that the maximum relative error is very small and it shows the ensemble has tracked the regulation signal at every time step (even after relaxing the constraint set to convex and rounding off later to find the device status). The distribution of relative error percentage corresponding to every regulation signal with and without additional cost on the decision variable $s$ is shown in Fig. \ref{fig:full_error_hist}. 

We next discuss the difference between the existing VB parameter identification approach such as \cite{hao2015aggregate} and draw comparisons with this work. The initial state of VB in \cite{hao2015aggregate} is considered to be zero and this assumption leads to a conservative estimate of the VB corresponding to the ensemble considered. Further in \cite{hughes2016identification}, the authors make the initial \textit{soc} an optimization variable which also leads to conservative VB parameters. To see the conservatism and dependency of VB parameters on initial condition, we consider two different scenarios for computation of VB parameters. The initial condition in the first scenario is considered to be zero and in the later one, the initial state of VB is computed analytically as shown in Eq. \eqref{eq:VB_AC_ic0} for AC devices and Eq. \eqref{eq:VB_WH_ic0} for EWH devices. 
\begin{table}[!ht]
\centering
\begin{tabular}{|c|c|c|c|c|c|}
\hline
\textbf{Type} & \textbf{\# ACs} & \textbf{\# EWHs} & $\textbf{a}$ (1/h) & $\textbf{C}$ (KWh)  & $x_{0}$ (KWh)   \\ \hline
Homogeneous   & 100 &  0  & 0.0328 & 56.181 &  0   \\ \hline
Homogeneous   & 100 &  0  & 0.033  & 61.515 &  43.682   \\ \hline
Homogeneous   &  0  & 120 &  0.0855    & 71.954 &  0   \\ \hline
Homogeneous   &  0  & 120 & 0.262  & 80.913 &   71.508  \\ \hline
Heterogeneous & 100 & 120 & 0.0987 & 126.869 & 117.807 \\ \hline
\end{tabular}
\caption{VB parameters for various ensembles}
\label{table:VB_parameters}
\end{table}

The initial condition from the analytic expression is zero only under a special case where the initial temperature of every device is equal to their lower (upper) point of the temperature hysteresis band of AC (EWH) devices. It can be seen that as expected, the zero initial condition for ACs and EWHs led to a smaller VB capacity as shown in Table \ref{table:VB_parameters}. The VB parameters shown in Table \ref{table:VB_parameters} are obtained by solving the optimization problem, Eq. \eqref{eq:VB_optm_problem} in Julia using Ipopt solver. Further, for computational convenience, it is assumed that $C_2 = -C_1 = C$. 

\begin{figure}
    \centering
    \includegraphics[width = 0.65 \linewidth]{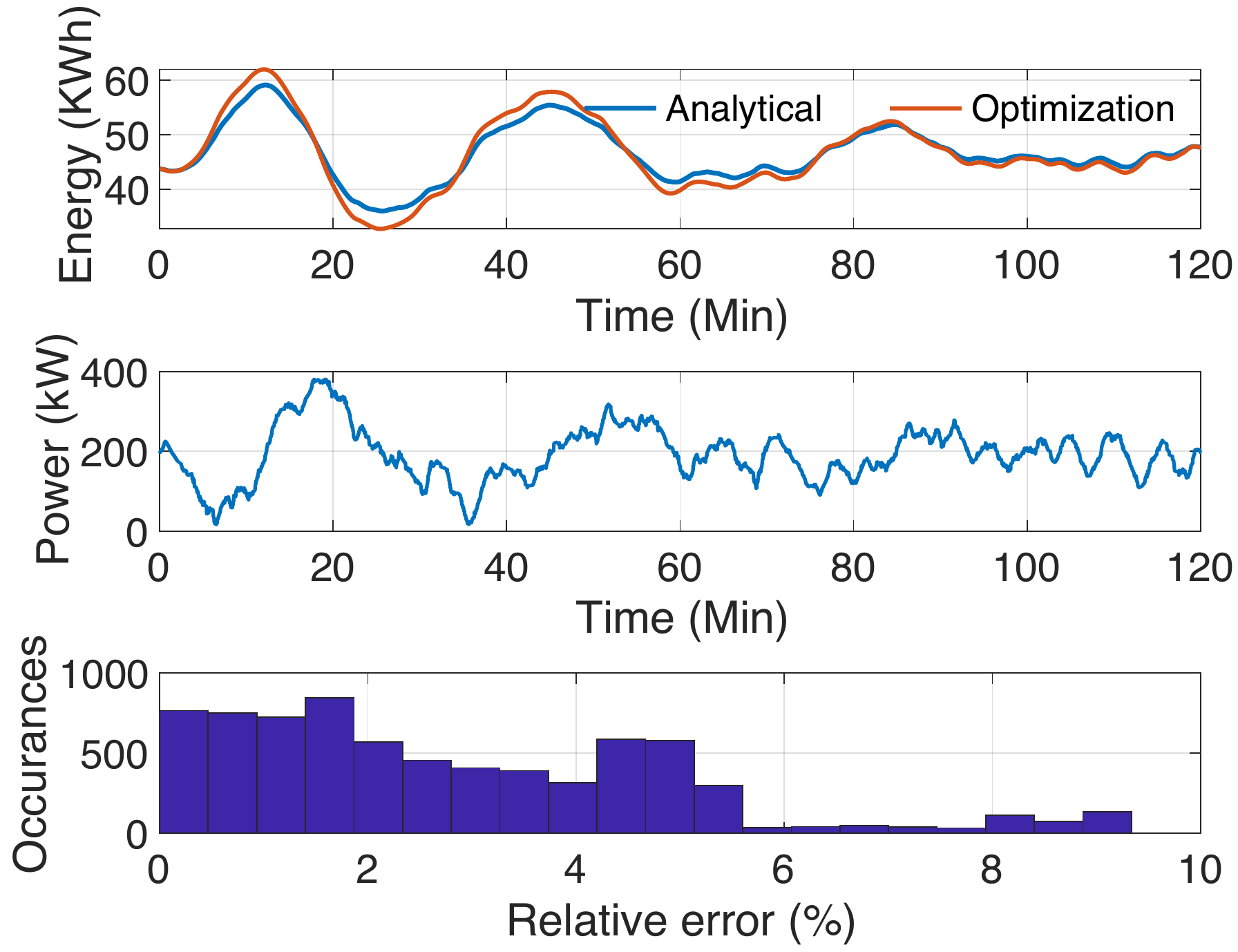}
    \caption{Virtual battery state comparison and RMS error}
    \label{fig:VB_state_comparison}
\end{figure}
We now present results that validate the VB parameters and we first consider the AC ensemble for this purpose. A regulation signal is considered and the VB state, energy corresponding to the AC ensemble is computed through the proposed framework. Also, the VB state is computed by extending Eq. \ref{eq:VB_AC_ic0} for all time points and we refer to this as the analytical VB state. 
Figure \ref{fig:VB_state_comparison} shows the comparison of the VB states computed using analytical and optimization based methods. Figure \ref{fig:VB_state_comparison} also show the regulation signal for which the VB state is computed and RMS error between the VB states obtained analytically and through the proposed approach.

\begin{figure}
    \centering
    \includegraphics[width = 0.65 \linewidth]{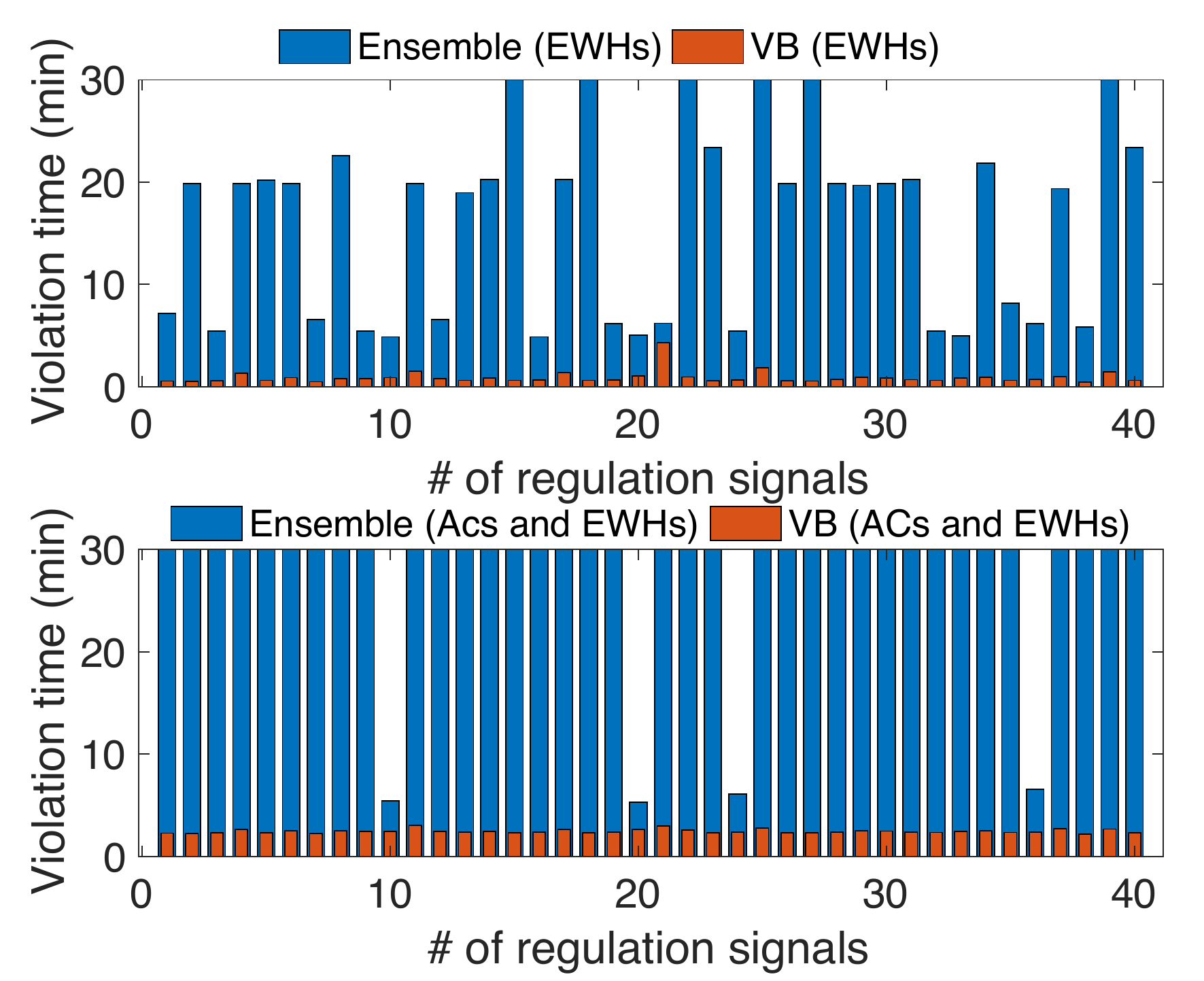}
    \caption{Violation times comparison}
    \label{fig:VB_verification}
\end{figure}
Next, we show a comparison of violation times between the actual ensemble following nonlinear dynamics and the first order linear VB system to further validate the VB models identified. Figure \ref{fig:VB_verification} shows the comparison of violation times for EWHs ensemble and also for the case of heterogeneous ensemble (ACs and EWHs). It can be seen that the violation times for the VB tend to be very small due to the fact that there are only two optimization variables ($a$ and $C$) for any given set of regulation signals and for every regulation signal, the violation time constraint (constraint in Eq. \eqref{eq:VB_optm_problem}) has to be satisfied. 

\subsection{VB Model: Heterogeneous Devices}
This section consists of VB parameter identification for heterogeneous devices. In this work, we restrict ourselves to a heterogeneous ensemble with AC and EWH devices. However, with less difficulty, it can be extended to a case with real battery, electric vehicle and solar PV generation along with ACs and EWHs. 

\begin{figure}
    \centering
    \includegraphics[width = 0.65 \linewidth]{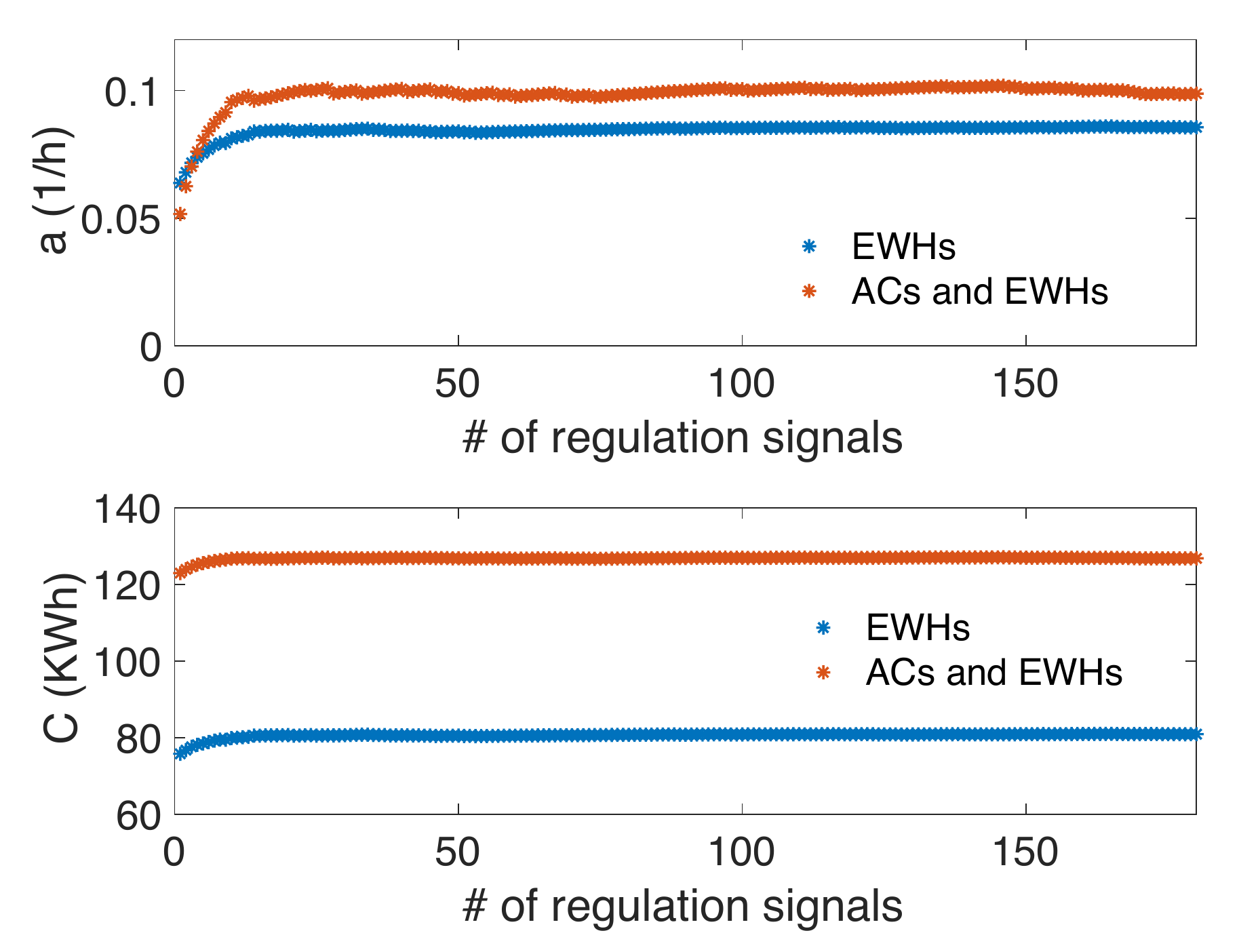}
    \caption{Evolution of VB parameters}
    \label{fig:VB_parameter_evolution}
\end{figure}
The heterogeneous ensemble consists of $100$ AC devices and $120$ EWH devices. The time-varying power limits are computed by applying the one-sided binary search algorithm discussed in Section \ref{sec:VB_model_identification}. The regulation signals which lie inside the power limits are chosen and violation times of the ensemble corresponding to those regulation signals are identified. The violation times are computed by solving the control design optimization problem, Eq. \eqref{eq:optim_problem_convex}. 

The violation times are then used to compute the VB parameters such as self dissipation rate and capacity. The initial state of the VB is identified by adding the initial state corresponding to the AC and EWH ensemble (see Eqs. \eqref{eq:VB_AC_ic0} and \eqref{eq:VB_WH_ic0}). Figure \ref{fig:VB_parameter_evolution} show that the VB parameters such as self dissipation and capacity saturate when the number of regulation signals are increased.

\section{Concluding Remarks and Future Scope}
\label{sec:conclusion}
This work consists of VB parameters computation for a heterogeneous ensemble of devices such as ACs and EWHs. A control design optimization problem is developed which identifies the device status (ON/OFF) in the ensemble while tracking the regulation signal and simultaneously respecting the device level temperature constraints. The initial state of the VB is identified analytically and time-varying power limits are computed using the proposed binary search algorithm. The self dissipation rate and capacity of the VB are obtained as a solution to the optimization problem which minimizes the tracking performance of the ensemble and first order VB system. Finally, by applying the proposed framework, VB parameters for  homogeneous/heterogeneous ensembles are identified and an extensive validation study shows the robustness and effectiveness of the proposed method. Further, it overcomes the conservatism with existing VB parameter identification methods. 

In the subsequent publications, we would like to extend the VB computation with ACs and EWHs to include real batteries, electric vehicles and solar PV inverters. Furthermore, the capacity bounds for the VB can be made more tighter by considering the lower and upper bounds of the VB capacity as decision variables. Moreover, finding the device parameters with some accuracy is always a challenge and remains to be explored. 

\bibliographystyle{unsrt}  
\bibliography{references}  






\end{document}